\documentclass[12pt]{article}
\usepackage{amsmath}
\usepackage{amssymb}
\usepackage[mathscr]{eucal}
\usepackage[usenames]{color}
\usepackage{url}
\usepackage[utf8]{inputenc}
\usepackage{graphicx}
\usepackage{epstopdf}
\usepackage{setspace}
\usepackage{color}
\usepackage{subfigure}
\usepackage{tikz}
\usepackage{authblk}

\definecolor{webgreen}{rgb}{0,0.5,0}
\newtheorem{rmk}{Remark}

\numberwithin{equation}{section}
\numberwithin{thm}{section}
\numberwithin{lemma}{section}
\numberwithin{prop}{section}
\numberwithin{cor}{section}
\numberwithin{rmk}{section}
\numberwithin{defn}{section}

\definecolor{darkolivegreen}{rgb}{0.333333, 0.419608, 0.1843140}

\begin{document}
%\pagenumbering{gobble}% Remove page numbers (and reset to 1)
\pagenumbering{arabic}

\title{\Large Nonlinear Opinion Dynamics in Systems with \\Leadership Effect}

\author[1,2]{F. Ata\c s\thanks{atasf@itu.edu.tr, fatma.atas@turkcell.com.tr}}
\author[1]{A. Demirci \thanks{demircial@itu.edu.tr}}
\author[1]{C. \"{O}zemir \thanks{ozemir@itu.edu.tr}}
\affil[1]{\small Department of Mathematics, Faculty of Science and Letters,\\ \small Istanbul Technical University, 34469 Istanbul, Turkey}
\affil[2]{Turkcell Technology Research and Development Company, Istanbul, Turkey }

\renewcommand\Authands{ and }

\date{\today}
\clearpage
\maketitle
\thispagestyle{empty}

\begin{abstract}
Opinion dynamics of a group of individuals is the change in the members' opinions through mutual interaction with each other. The related literature contains works in which the dynamics is modeled as a continuous system, of which behavioral patterns are analyzed in regard to the parameters contained in the system. These models are constructed by the assumption that the individuals are independent. Besides, the decisions of the individuals are only affected by two forces: self-bias force and group influence force. In this work we consider a nonlinear dynamical system which models the evolution of the decision of a group under the existence of a leader. Bifurcation analysis of the system is performed to obtain stability results on the system.
\end{abstract}

\section{Introduction}

Among the fields that the theory of dynamical systems finds its applications one can count the  models of swarm behaviour and consensus processes. Problems of this kind appear in physics of many natural phenomena and engineering applications: Immigrating groups of animals (flocks of birds and fish), robots that are in communication with each other on a bounded region, motions of groups of stars or galaxies, diffusive dynamics of microorganisms. The basic rule in the dynamics of the group consisting of $N$  individuals is that each agent moves in a spacetime path such that they are far enough from each other so they do not collide and keep a distance close enough so as to keep in communication with the neighboring agents, which corresponds to an evolution under an attractive-repulsive potential. In this setting, one of the simplest models is the ODE system
\begin{equation}
\frac{dx_i}{dt}=\sum_{j=1}^N a_{ij}(x_j-x_i), \qquad  i=1,2,...,N.
\end{equation}
In case of collective motion, $x_i$ stands for the position of an agent in 1-dimension. In opinion dynamics, which is the  evolution of opinions of individuals in a society through mutual interaction with other individuals, $x_i$ is the value of the "opinion" that the agent $i$ produces.

One of the questions to be answered in consensus problems is whether the agents in communication will arrive at a common decision
in the limit $t\rightarrow \infty$ or not. Or, are there subgroups with different opinion values in the long run? There are many variants of this model in literature \cite{Saber2007,Motsch2014,Bertozzi2012}. For a substantial review on applications of the theory of dynamical systems to decision dynamics, models of cultural change, language evolution, biological evolution and swarm motion one can see \cite{Castellano2009}.

The model in which the opinions of $N$ individuals are affected by their initial self-judgements and the weighted average of the group's opinion in discrete time is known as Friedkin- Johnsen model \cite{Friedkin2011}:
\begin{equation}\label{FJ}
x_i (k+1)=a_i \sum_{j=1}^N w_{ij} x_j (k)+(1-a_i ) x_i (0);   \quad     i=1,...,N,   \quad   k\geq 0.
\end{equation}
Here $x_i (k)$ is the opinion of the $i$th individual at time $k$,  $x_i (0)$ is its opinion at the beginning, $a_i$ is its sensitivity and $w_{ij}$ is a measure of the influence of agent $j$ on the decision of agent $i$. If one subtracts the term  $x_i (k)=(1+a_i-a_i)x_i (k)$ from both sides of \eqref{FJ} and organize the terms,
\begin{equation}\label{discrete}
 x_i (k+1)-x_i (k)=a_i \sum_{j=1}^N w_{ij} (x_j (k)-x_i (k))-(1-a_i )(x_i (k)-x_i (0))
\end{equation}
is  obtained \cite{Gabbay2014}.

In \cite{Gabbay2014}, instead of the discrete model \eqref{discrete}, the continuous counterpart
\begin{equation}\label{gab}
\frac{dx_i}{dt}=-\gamma_i (x_i-\mu_i )+\sum_{j=1}^N \kappa_{ij} h(x_j-x_i ), \quad         i=1,...,N
\end{equation}
is analyzed for long-time behavior according to the parameters included. Here the opinions of the agents in a group change under the influence of two main forces: One is the self-bias force, the other is the group influence force.
Self-bias force represents the agent's dependency on its initial opinion and is taken into account by the first term on the righthand side of \eqref{gab}. With $\gamma_i>0$, this force has a negative impact on the change of opinion. $\gamma_i$ is the strength of the $i$th individual's dependency on its initial opinion $\mu_i$. Group influence force is a result of the effect of other individuals on a specific agent, being represented by the second term on the righthand side of \eqref{gab}. $\kappa_{ij}>0$ is the coupling strength and is proportional to the frequency in which agent $j$ communicates with agent $i$.

$h(x_j-x_i )$  is the coupling function and depending on the opinion difference $x_j-x_i$ of the agents $i$ and  $j$, represents the measure of the strength of influence that agent $j$ has on agent   $i$. Due to the term $h(x_j-x_i )$  the model is a nonlinear one, taken in the form
\begin{equation}
h(x_j-x_i )= (x_j-x_i )\exp\Big[-\frac{1}{2}  \frac{(x_j-x_i )^2}{\lambda_i^2 }\Big].
\end{equation}
Therefore, the effects of individuals on each other, of which opinions are far away from each other,  are taken to be as less and less.

\cite{Gabbay2014} considers  the system \eqref{gab}  for a group of $N=3$ agents, explicitly in the form
\begin{align}
&\frac{dx_1}{dt}=-(x_1-\mu_1)+(\kappa+\nu)h(x_2-x_1),\nonumber\\
&\frac{dx_2}{dt}=-x_2+(\kappa-\nu)(h(x_1-x_2)+h(x_3-x_2)), \label{21}\\
&\frac{dx_3}{dt}=-(x_3-\mu_3)+(\kappa+\nu)h(x_2-x_3).\nonumber
\end{align}

\begin{figure}[!ht] \centering
\subfigure[$\kappa=1$]{\label{sek11} %% label for first subfigure
\includegraphics[scale=0.50]{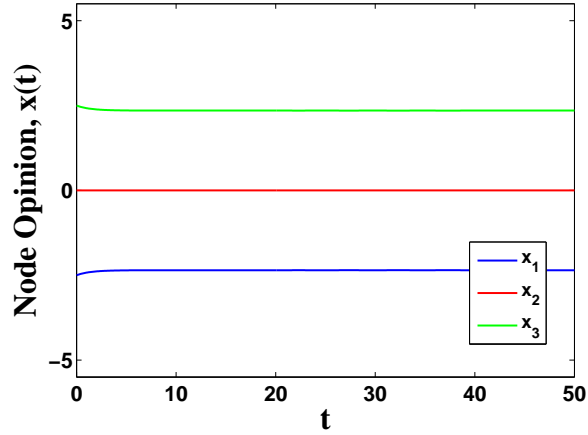}} \hspace{2cm}
\subfigure[$\kappa=1.5$]{
\label{sek12} %% label for second subfigure
\includegraphics[scale=0.5]{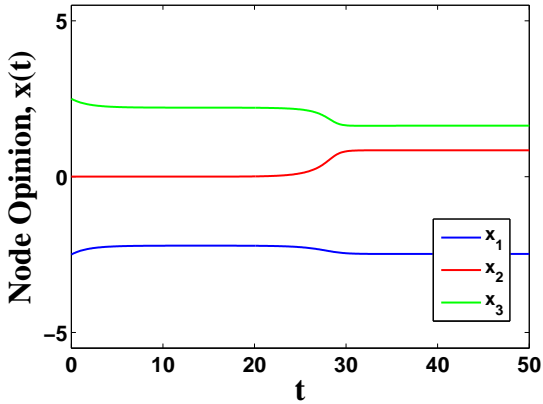}}
\subfigure[$\kappa=3$]{\label{sek13} %% label for first subfigure
\includegraphics[scale=0.5]{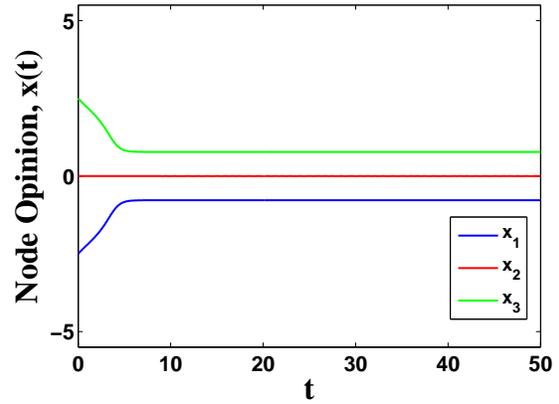}}\hspace{2cm}
\caption{{Equilibrium outcomes in the dynamical system \eqref{21}  with  initial conditions: $x_1(0)=-2.5$, $x_2(0)=10^{-6}$, $x_3(0)=2.5$ and  initial disagreement: $\Delta\mu=\mu_3-\mu_1=5$. }}
\label{gabb} %% label for entire figure
\end{figure}

Figure \ref{gabb}, a reproduction of the findings in \cite{Gabbay2014} on numerical solution of \eqref{21},  indicates  the change that occurs on the opinion of agents in a chain network for three different values of the  coupling strength $\kappa$. They observe three distinct states which correspond to these values of the coupling strength. For $\kappa=1$, Figure \ref{sek11} demonstrates a state of Symmetric High Discord (SHD). In this state, the center node stays at zero and the natural biases of end nodes change only a small amount. For $\kappa=1.5$, Figure \ref{sek12} indicates Majority Rule state. In this state, center node approaches one of the end nodes generating a majority rule pair.
For $\kappa=3$, Figure \ref{sek13} shows the state of Symmetric Low Discord. In this state, the end nodes approach to the decision of the center node. Let us note that for  \cite{Gabbay2014} employs an initial condition of $x_2(0)=10^{-6}$ instead of a zero initial condition so as to create a small perturbation for $x_2$ in the direction of $x_3$.

In SHD state, the discord of members in the group is apparent. It is impossible that all group members agree to the same decision. In  MR state, nodes that form a majority pair accept a mutual decision. Although majority rule is adequate to reach a decision, it is not the decision of the group. If majority rule is sufficient for arriving at a decision, it will be the \emph{policy of the group} \cite{Gabbay2014}.  In SLD state, all group members agree on a common decision. This decision is very close to the opinion of center node, albeit we do not give any precise definition of this "closeness".

The long-time behaviour (the stable picture) of the system is analyzed in the parameter space $(\Delta \mu,\kappa)$; where $\Delta \mu=\mu_3-\mu_1$ is the maximum initial opinion difference between agents and $\kappa$ is a uniform coupling strength, to distinguish the behaviours (a) symmetric high discord-SHD,  (b) majority rule-MR, (c) symmetric low discord-SLD. The information given by the analysis according to the parameters included in the system, without solving the system, is the final consensus level of the group: Disagreement (or \emph{deadlock}  \cite{Gabbay2014}), majority, or a level acceptable as an agreement.

Another occasion of interest to decision dynamics is the situation in which the group under consideration is including a leader. In such systems, the leader is defined as the power that controls the opinion dynamics of the group. Hegselmann- Krause (HK) model is an example to these systems \cite{HK2014}:
\begin{subequations}\label{HK}
\begin{align}
&\frac{dx_0}{dt}=u(t),\\
&\frac{dx_i}{dt}=\sum_{j=1}^N a_{ij}(x_j-x_i)+C_i(x_0-x_i), \quad i=1,2,...,N.
\end{align}
\end{subequations}
Here the leader is a "stubborn" agent which is not influenced by anyone but has an influence on everyone. As it is seen, the individual $x_0$ has a dynamics determined by the function $u(t)$, not affected by any other agent, and has an effect in each agent's steering equation and represents the leader of the group \cite{HK2014}.

With the motivation provided by the literature outlined above, this article is focused on the problem framed by the following questions:
\begin{itemize}
\item How can the opinion evolution given in \eqref{gab} be modeled in the presence of a leader? This is an open decision dynamics problem not analyzed in the existing literature, to the best of our knowledge, and the answer starts by considering systems \eqref{gab} and \eqref{HK} together.
\item Through construction of such a model, and by a similar machinery as in \cite{Gabbay2014} (via bifurcation analysis in parameter space), can we foresee and identify the final stable situation of this system as disagreement, majority or agreement?
\end{itemize}
When the systems \eqref{gab} and \eqref{HK} are considered together, the model that will be the base of our problem is of the form
\begin{subequations}\label{sysfull}
\begin{align}
&\frac{dx_0}{dt}=u(t),\\
&\frac{dx_i}{dt}=-\gamma_i(x_i-\mu_i)+\sum_{j=1}^N\kappa_{ij}h(x_j-x_i)+C_i(x_0-x_i), \quad i=1,2,...,N;
\end{align}
\end{subequations}
which is a combination of the dynamical systems considered in \cite{Gabbay2014} and \cite{HK2014}. Just to mention, $x_0(t)$ is the opinion function of a stubborn agent, called the  leader, as it is not affected by the opinion of anyone in the group while he/she influences the opinion of everyone in the group. The other individuals have the same dynamics in the system \eqref{gab} besides  being influenced from the evolution of $x_0$.  $C_i$ represents the strength of this influence on each agent \cite{HK2014}. This constant is utilized to tune the free decision maker's, the leader's effect on the individuals of the society.
\section{Analysis}
For simplicity we will consider only the case $u(t)\equiv 0$, which gives a leader at a constant decision $x_0$ at any time. Therefore, we are led to the dynamical system
\begin{equation}\label{sys}
\frac{dx_i}{dt}=-\gamma_i(x_i-\mu_i)+\sum_{j=1}^N\kappa_{ij}h(x_j-x_i)+C_i(x_0-x_i),
\end{equation}
for $i=1,2,...,N$. The last term in \eqref{sys} represents the leadership effect. See that if $x_0>x_i$ this term supports $x_i$ to increase to arrive at $x_0$, and vice versa if $x_0<x_i$, when  $C_i>0$. We see that this term works for making $x_i$ come and stay around the opinion $x_0$. In case $C_i<0$, the effect is to push $x_i$  away from the opinion $x_0$. The term with $C_i$ can also be viewed as a control term and $C_i$ as the control parameter. In what follows we will always consider $C_i$ as a  positive parameter to represent the dominance of the leader on the group and impose it on the system as $\pm C_i$ for an attractive or impulsive leader effect.

Following \cite{Gabbay2014}, we analyze the evolution of opinion on a network which has three nodes. The natural biases in this network are symmetric around zero: $\mu_1=-\Delta\mu/2$, $\mu_2=0$ and $\mu_3=\Delta\mu/2$. The network has a chain topology; that is, node 1 and 3 are not conjoined while both node 1, node 2 and node 2, node 3 are connected between each other. Therefore, node 1 and node 3 do not directly influence each other. The binary adjacency matrix is then
\begin{equation}A=
\begin{bmatrix}
0&1&0\\1&0&1\\0&1&0
\end{bmatrix}.
\end{equation}
Based on this topology, \eqref{sys} takes the form
\begin{align}\label{sys1}
&\frac{dx_1}{dt}=-(x_1-\mu_1)+(\kappa+\nu)h(x_2-x_1)+C_1(x_0-x_1),\nonumber\\
&\frac{dx_2}{dt}=-x_2+(\kappa-\nu)(h(x_1-x_2)+h(x_3-x_2))+C_2(x_0-x_2)\\
&\frac{dx_3}{dt}=-(x_3-\mu_3)+(\kappa+\nu)h(x_2-x_3)+C_3(x_0-x_3).\nonumber
\end{align}
In this setting, asymmetric coupling between the center node 2 and the end nodes is possible by virtue of the parameter $\nu$. If $\nu$ is a positive parameter then  node 2 has a bigger effect on the end nodes. If $\nu$ is negative the effect of the  end nodes is higher.

We define three new parameters $r$, $s$ and $\bar{x}$: $r$ is described as discord and expressed as $r=x_3-x_1$. $s$ is asymmetry of the opinions, defined as  $s=x_3-x_2-(x_2-x_1)=x_3-2x_2+x_1$. $\bar{x}$ is the average of opinions (mean node opinion): $\bar{x}=\frac{1}{3}(x_1+x_2+x_3)$. The  coupling function $h(x)$ is an odd function. In terms of these variables, \eqref{sys1} takes the form
\begin{align}
\frac{dr}{dt}&=-r+\mu_3-\mu_1+(C_3-C_1)x_0+C_1x_1-C_3x_3-(\kappa+\nu)(h(\frac{r+s}{2})+h(\frac{r-s}{2})),\nonumber\\
\frac{ds}{dt}&=-s+(C_3-2C_2+C_1)x_0-C_1x_1+2C_2x_2-C_3x_3+\mu_3+\mu_1\nonumber\\
              &-(3\kappa-\nu)(h(\frac{r+s}{2})-h(\frac{r-s}{2})), \label{rsx}\\
\frac{d\bar{x}}{dt}&=-\bar{x}+\frac{1}{3}(C_1+C_2+C_3)x_0+\frac{1}{3}(\mu_3+\mu_1)
                    -\frac{1}{3}(C_1x_1+C_2x_2+C_3x_3)\nonumber\\
                    &-\frac{2\nu}{3}(h(\frac{r+s}{2})-h(\frac{r-s}{2})).\nonumber
\end{align}

\begin{figure}[!ht] \centering
\subfigure[$\kappa=0.5$]{\label{sek21} %% label for first subfigure
\includegraphics[scale=0.50]{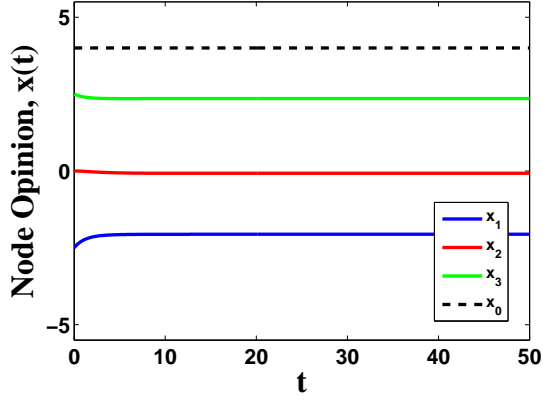}} \hspace{2cm}
\subfigure[$\kappa=1.5$]{\label{sek22} %% label for second subfigure
\includegraphics[scale=0.50]{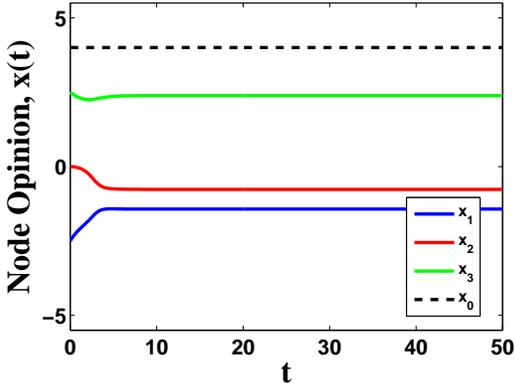}}
\subfigure[$\kappa=14$]{\label{sek23} %% label for first subfigure
\includegraphics[scale=0.56]{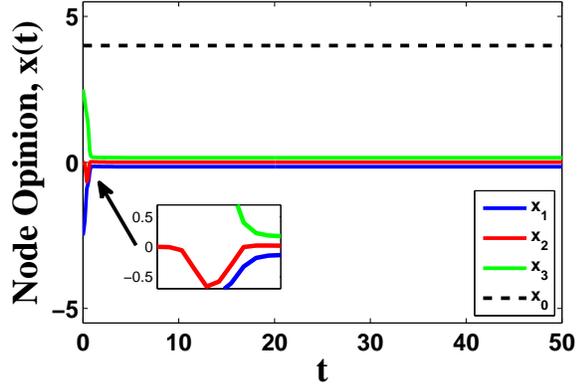}}\hspace{2cm}
 \caption{{Equilibrium outcomes for the system \eqref{main} with control/leadership effect. $x_0=4$, $C=0.05$. Initial disagreement: $\Delta\mu=5$ and initial conditions: $x_1(0)=-2.5$, $x_2(0)=0$, $x_3(0)=2.5$.  }}
\label{imper} %% label for entire figure
\end{figure}

We shall assume that the leader applies a pulling force on node 1 while it applies a repelling force on node 3. That is, we assume that $C_1=C$, $C_2=0$ and $C_3=-C$ where $C>0$. The equations of motion in the variables $x_1$, $x_2$, $x_3$, just to stress, are
\begin{equation}\label{main}
\begin{split}
\frac{dx_1}{dt}&=-(x_1-\mu_1)+(\kappa+\nu)h(x_2-x_1)+C(x_0-x_1) \\
\frac{dx_2}{dt}&=-x_2+(\kappa-\nu)(h(x_1-x_2)+h(x_3-x_2))\\
\frac{dx_3}{dt}&=-(x_3-\mu_3)+(\kappa+\nu)h(x_2-x_3)-C(x_0-x_3)
\end{split}
\end{equation}
Figure \ref{imper} shows  the change of opinion by using three distinct values of the coupling strength $\kappa$. SHD, MR and SLD states correspond to $\kappa=0.5$, $\kappa=1.5$ and $\kappa=14$, respectively.

From these simulations we see that it is also possible to observe the the three different states of SHD, MR and SLD of \cite{Gabbay2014} for the modified system \eqref{main} with the specific leadership effect considered in our work.
In SHD state, the initial opinions of nodes in the group do not change almost at all, since the value of interaction forces, coupling strength and leadership effect, are very small. In MR state, the center node generates majority pair by approaching to one of end nodes. It selects node 1 to generate the majority pair because of the leadership effect. In SLD state, the distance  between the center node $x_2$ and end nodes $x_1$ and $x_3$ decreases a considerable amount.

When we rearrange \eqref{rsx}  for $C_1=C$, $C_2=0$ and $C_3=-C$ we obtain
 \begin{subequations}
\begin{align}
\label{33a}\frac{dr}{dt}&=-r+2C\bar{x}+\frac{C}{3}s-2C x_0+\mu_3-\mu_1-(\kappa+\nu)(h(\frac{r+s}{2})+h(\frac{r-s}{2})), \\
\label{33b}\frac{ds}{dt}&=-s+Cr+\mu_3+\mu_1-(3\kappa-\nu)(h(\frac{r+s}{2})-h(\frac{r-s}{2})),\\
\label{33c}\frac{d\bar{x}}{dt}&=-\bar{x}+\frac{C}{3}r+\frac{1}{3}(\mu_3+\mu_1)-\frac{2\nu}{3}(h(\frac{r+s}{2})-h(\frac{r-s}{2})).
\end{align}
\end{subequations}
where we have used the fact that $\displaystyle x_1+x_3=2\bar{x}+s/3$.

In the sequel, bifurcation analysis will be made in the chain network. The system will be investigated in regard to states identified above and the passage from one  state to another will analyzed. The most interesting of all, the existence of an imperfect pitchfork bifurcation will be shown and several bifurcation boundaries will be identified. While the leaderless system exhibits a pitchfork bifurcation, the perturbation of the model by a leader causes the qualitative behaviour to be perturbed to the imperfect one.

\eqref{33a} is examined in the equilibrium of SHD state ($s=0$). It means that the system does not have asymmetric coupling in this state. Furthermore, mean node opinion is zero ($\bar{x}=0$) in this state and the value of discord is accepted as $r\approx\Delta\mu+\theta$ for large $\Delta\mu$. \eqref{33a} is stated for $\bar{x}=0$ and $s=0$ as below
\begin{equation}
\label{34}\frac{dr}{dt}=-r-2C x_0+\mu_3-\mu_1-2(\kappa+\nu)h(\frac{r}{2}).
\end{equation}
We rearrange this equation for $r\approx\Delta\mu+\theta$ as
\begin{equation}
\label{35}0=\theta+2C x_0+2(\kappa+\nu)h(\frac{\Delta\mu+\theta}{2}).
\end{equation}
The Taylor expansion of the coupling function is found around $\theta=0$ in order to calculate the value of $\theta$ as follows
\begin{equation}
h(\frac{\Delta\mu+\theta}{2})\approx h(\frac{\Delta\mu}{2})+h'(\frac{\Delta\mu}{2})\frac{\theta}{2}.
\end{equation}
\eqref{35} is rearranged by utilizing this expansion to find
$$\theta=-\frac{2Cx_0+2(\kappa+\nu)h(\frac{\Delta\mu}{2})}{1+(\kappa+\nu)h'(\frac{\Delta\mu}{2})}$$
where $\displaystyle h(\frac{\Delta\mu}{2})=\frac{\Delta\mu}{2}e^{-\frac{\Delta\mu^{2}}{8}}$ and
$\displaystyle h'(\frac{\Delta\mu}{2})=\frac{1}{2}(1-\frac{\Delta\mu^{2}}{4})\ e^{-\frac{\Delta\mu^{2}}{8}}$.
In order to exhibit the existence of bifurcation, we consider small perturbations of s around $s=0$ in \eqref{33b}. The Taylor approximation of $h(\frac{r+s}{2})-h(\frac{r-s}{2})$ is
\begin{equation}h(\frac{r+s}{2})-h(\frac{r-s}{2})\approx h'(\frac{r}{2})s+\frac{1}{24}h'''(\frac{r}{2})s^{3}.\end{equation}
We rewrite \eqref{33b} by using this expansion as below
\begin{equation}
\label{38}\frac{ds}{dt}\approx Cr-(1+(3\kappa-\nu)h'(\frac{r}{2}))s-\frac{1}{24}(3\kappa-\nu)h'''(\frac{r}{2})s^{3}.
\end{equation}
\eqref{38} is rescaled with the parameters $\tau$ and $R$ as
\begin{subequations}
\begin{align}
\tau=[\frac{1}{24}(3\kappa-\nu)&h'''(\frac{r}{2})]t, \quad R=-\frac{1+(3\kappa-\nu)h'(\frac{r}{2})}{\frac{1}{24}(3\kappa-\nu)h'''(\frac{r}{2})},\\
\label{39b}&\frac{ds}{d\tau}=A+Rs-s^{3}
\end{align}
\end{subequations}
where $\displaystyle A=\frac{Cr}{\frac{1}{24}(3\kappa-\nu)h'''(\frac{r}{2})}$. \eqref{39b} is nothing but the general form of an imperfect pitchfork bifurcation.  The critical point of \eqref{39b} will be obtained for $\Delta=0$, where $\Delta$ is the discriminant of the cubic polynomial on the righthandside of \eqref{39b}. This condition is going to give us one of the boundaries separating the different regions in the stability diagram Figure 3.
\begin{rmk}
The main outcome that makes this paper a contribution to the literature is the emergence of  Eq. \eqref{39b}. In the leaderless case, it had been reported in \cite{Gabbay2014} that the dynamical model exhibits a pitchfork bifurcation with the prototype equation for $s$ being \eqref{39b} with $A=0$. In our case, the prototype equation for $s$ becomes the more general one \eqref{39b} and the system undergoes an imperfect pitchfork bifurcation in the existence of a leader.
\end{rmk}
\subsection{Boundary $\kappa_1$}
$\kappa_1$, depicted in Figure 3,  is the upper boundary of SHD in $\Delta\mu-\kappa$ parameter space. This boundary is found by using the critical point of the imperfect pitchfork bifurcation. This critical case is the point  where the righthand side of \eqref{39b} has a unique real root before and has three real roots afterwards, or vice versa. To mention,   for the general cubic equation  $ax^{3}+bx^{2}+c x+d=0$ the discriminant $\Delta$ is expressed as follows
\begin{equation}
\label{310}\Delta=b^{2}c^{2}-4ac^{3}-4b^{3}d-27a^{2}d^{2}+18abcd.
\end{equation}
It is clear that the critical point of the imperfect bifurcation is obtained when $\Delta=0$. To obtain $a, b, c, d$, \eqref{38} is rearranged for $r\approx\Delta\mu+\theta$ as
\begin{equation}
\frac{ds}{dt}\approx C(\Delta\mu+\theta)-(1+(3\kappa-\nu)h'(\frac{\Delta\mu+\theta}{2}))s-\frac{1}{24}(3\kappa-\nu)h'''(\frac{\Delta\mu+\theta}{2})s^{3}.
\end{equation}
By using the coefficients of above equation, $a, b, c, d$ are obtained as below
\begin{equation}
\begin{split}
&a=-\frac{1}{24}(3\kappa-\nu)h'''(\frac{\Delta\mu+\theta}{2}),\,\,\,\,\,\,\,\,\,\,b=0,\\
&c=-(1+(3\kappa-\nu)h'(\frac{\Delta\mu+\theta}{2})),\,\,\,\,\,d=C(\Delta\mu+\theta)
\end{split}
\end{equation}
These parameters are substituted in $\Delta = 0 $ to get
\begin{equation}\label{313c}
-32(1+(3\kappa-\nu)h'(\frac{\Delta\mu+\theta}{2}))^{3}=9C^{2}(\Delta\mu+\theta)^{2}(3\kappa-\nu)h'''(\frac{\Delta\mu+\theta}{2}).
\end{equation}
The Taylor approximations of $h'(\frac{\Delta\mu+\theta}{2})$ and $h'''(\frac{\Delta\mu+\theta}{2})$ are calculated around $\theta=0$ as follows
\begin{subequations}
\begin{align}
&h'(\frac{\Delta\mu+\theta}{2})\approx h'(\frac{\Delta\mu}{2})+h''(\frac{\Delta\mu}{2})\frac{\theta}{2},\\
&h'''(\frac{\Delta\mu+\theta}{2})\approx h'''(\frac{\Delta\mu}{2})+h^{(4)}(\frac{\Delta\mu}{2})\frac{\theta}{2}.
\end{align}
\end{subequations}
These expansions are substituted in \eqref{313c} as
\begin{align}
\label{315}-32(1+(3\kappa-\nu)\Big[ h'(\frac{\Delta\mu}{2})+h''(\frac{\Delta\mu}{2})\frac{\theta}{2}\Big]^{3}=9C^{2}(\Delta\mu+\theta)^{2}(3\kappa-\nu)\Big[ h'''(\frac{\Delta\mu}{2})+h^{(4)}(\frac{\Delta\mu}{2})\frac{\theta}{2}\Big]
\end{align}
where
\begin{subequations}\label{219}
\begin{align}
h'''(\frac{\Delta\mu}{2})   &=\frac{1}{16}(-6+3\Delta\mu^{2}-\frac{\Delta\mu^{4}}{8})e^{-\frac{\Delta\mu^{2}}{8}},\\
h^{(4)}(\frac{\Delta\mu}{2})&=\frac{1}{64}(30\Delta\mu-5\Delta\mu^{3}+\frac{\Delta\mu^{5}}{8})\ e^{-\frac{\Delta\mu^{2}}{8}}.
\end{align}
\end{subequations}
$\kappa_1$ is given in Figure 3 of the next Section, by the numerical solution of  \eqref{315} for $\kappa$ with varying values of $\Delta \mu$, making use of \eqref{219}.
\subsection{Boundary $\kappa_2$}
Now, following \cite{Gabbay2014}, we are going to analyze the situation where the system passes from MR state to SHD state, as the coupling constant is decreased. In the final picture Figure 3, this boundary $\kappa_2$ will be the curve above which MR and SHD both occur and below which the system is in SHD. The transition mechanism is a  saddle-node bifurcation, and $\kappa_2$ is a perturbed form of  that found in \cite{Gabbay2014}.
Stable and unstable equilibrium points collide with one another along this curve. To find this boundary, \eqref{33b} is analyzed around MR equilibrium point. In MR state, nodes $x_1$ and $x_2$  create a majority pair and $x_3$ is the minority node. $\displaystyle \varepsilon=\frac{x_1+x_2}{2}+\frac{\Delta\mu}{4}$ is defined as the change in the mean opinion of $x_1$ and $x_2$, and for large $\Delta\mu$, the term $\displaystyle h(\frac{r+s}{2})$ is negligible. The approximate values of $x_1$ and $x_3$  are assumed to be
\begin{equation}
x_1\approx-\frac{\Delta\mu}{2}+2\varepsilon-x_2, \qquad  x_3\approx\frac{\Delta\mu}{2}.
\end{equation}
The asymmetry is
\begin{equation}
s=x_3-2x_2+x_1=\frac{\Delta\mu}{2}-2x_2-\frac{\Delta\mu}{2}+2\varepsilon-x_2=-3x_2+2\varepsilon,
\end{equation}
which also gives
\begin{equation}
x_2=-\frac{s}{3}+\frac{2\varepsilon}{3}, \qquad  x_1=-\frac{\Delta\mu}{2}+\frac{s}{3}+\frac{4\varepsilon}{3}.
\end{equation}
The discord is
\begin{equation}\label{rr}
r=x_3-x_1=\frac{\Delta\mu}{2}+\frac{\Delta\mu}{2}-2\varepsilon-\frac{s}{3}+\frac{2\varepsilon}{3}=\Delta\mu-\frac{s}{3}-\frac{4\varepsilon}{3},\\
\end{equation}
yielding
\begin{equation}
\frac{r-s}{2}=\frac{1}{2}(\Delta\mu-\frac{s}{3}-\frac{4\varepsilon}{3}-s)=\frac{2}{3}(\frac{3}{4}\Delta\mu-s-\varepsilon).
\end{equation}
We define a new parameter as $\tilde{s}=\frac{3}{4}\Delta\mu-s-\varepsilon$ and when \eqref{33b} is rewritten for $\tilde{s}$ we get
\begin{equation}
\label{320}\frac{d\tilde{s}}{dt}=-(\tilde{s}-\frac{3}{4}\Delta\mu+\varepsilon)+Cr-(3\kappa-\nu)h(\frac{2}{3}\tilde{s}).
\end{equation}
\eqref{33c} in the equilibrium case (r.h.s. equal to zero) gives  the value of $\varepsilon$ as
\begin{equation}
\varepsilon=\frac{Cr}{2}+\nu h(\frac{2}{3}\tilde{s}).
\end{equation}
\eqref{320} is rearranged as below
\begin{equation}\label{322c}
\frac{d\tilde{s}}{dt}=-(\tilde{s}-\frac{3}{4}\Delta\mu)+\frac{C}{2}r-3\kappa h(\frac{2}{3}\tilde{s}).
\end{equation}
Considering  \eqref{33a} in the equilibrium case, we obtain
\begin{equation}\label{323d}
r(1-\frac{C}{2})=(C-1)\varepsilon+(1+\frac{C}{4})\Delta\mu-\frac{C}{3}\tilde{s}-2C x_0-\kappa h(\frac{2}{3}\tilde{s}).
\end{equation}
The value of $r$ in \eqref{rr} is reexpressed using $\tilde s$ as 
\begin{equation}
\label{324}r(1-\frac{C}{2})=\frac{3}{4}(1-\frac{C}{2})\Delta\mu+\frac{1}{3}(1-\frac{C}{2}){\tilde{s}}-(1-\frac{C}{2})\varepsilon.
\end{equation}
\eqref{323d} and \eqref{324} together gives
\begin{equation}
\varepsilon=-(\frac{1}{2C}+\frac{5}{4})\Delta\mu+\frac{1}{3}(1+\frac{2}{C})\tilde{s}+4 x_0+\frac{2}{C}\kappa h(\frac{2}{3}\tilde{s})].
\end{equation}
\\
\eqref{322c} is rewritten for $\displaystyle r=\frac{\tilde{s}}{3}+\frac{3}{4}\Delta\mu-\varepsilon$ to get
\begin{equation}\label{326b}
\frac{d\tilde{s}}{dt}=-(1-\frac{C}{6})\tilde{s}+\frac{3}{4}(1+C)\Delta\mu-\frac{C}{2}\varepsilon-3\kappa h(\frac{2}{3}\tilde{s}).
\end{equation}
The value of $\varepsilon$ is substituted in \eqref{326b} to obtain
\begin{equation}
\label{327}\frac{d\tilde{s}}{dt}=-\frac{4}{3}\tilde{s}+(1+\frac{11C}{3})\Delta\mu-2Cx_0-4\kappa h(\frac{2}{3}\tilde{s}).
\end{equation}
To reveal the saddle-node bifurcation in the equilibrium point, the right side of \eqref{327} and its derivative are set equal to zero. Thus,
\begin{subequations}
\begin{align}
\label{328a}&\frac{4}{3}\tilde{s}-(1+\frac{11C}{3})\Delta\mu+2Cx_0=-4\kappa_2\tilde{s}e^{-\frac{2}{9}\tilde{s}^{2}},\\
\label{328b}&\frac{4}{3}=-4\kappa_2(1-\frac{4}{9}\tilde{s}^2)e^{-\frac{2}{9}\tilde{s}^{2}}.
\end{align}
\end{subequations}
In order to obtain $\tilde{s}$, \eqref{328a} and \eqref{328b} are solved together and we get
\begin{equation}
\tilde{s}^{3}-\frac{3}{4}(\Delta\mu+\frac{11C}{8}\Delta\mu-2Cx_0)\tilde{s}^{2}+(\frac{27}{16}+\frac{297C}{128})\Delta\mu-\frac{54C}{16}x_0=0.
\end{equation}
$\tilde{s}=\frac{3}{2}+\frac{12}{(8+11C)\Delta\mu}$ is a solution of above equation. \eqref{328a} is solved with respect to $\kappa_2$ as follows
\begin{equation}
\label{330}\kappa_2=\frac{1}{4\tilde{s}}(\Delta\mu+\frac{11C}{8}\Delta\mu-\frac{4}{3}\tilde{s}-2Cx_0)e^{\frac{2}{9}\tilde{s}^{2}}
\end{equation}
The value of $\tilde{s}$ is substituted in \eqref{330} to obtain
\begin{align}
\kappa_2&=\frac{(8+11C)\Delta\mu}{6[(8+11C)\Delta\mu+8]}\Big(\Delta\mu+\frac{11C}{8}\Delta\mu-\frac{16}{(8+11C)\Delta\mu}-2Cx_0-2\Big)e^{\frac{2}{9}(\frac{3}{2}+\frac{12}{(8+11C)\Delta\mu})^{2}}. \end{align}

\subsection{Boundary $\kappa_3$}
Let $\kappa_3$ be the lower boundary of SLD state. To find this boundary, we will consider \eqref{34} in the equilibrium case and the coefficient of linear term in the \eqref{38}. A system which consists of these equations indicates the existence of a pitchfork bifurcation between SLD and MR states. This system is constructed as follows.
\begin{subequations}
\begin{align}
\label{332a}&r+2C x_0-\Delta\mu=-(\kappa_3+\nu)re^{-\frac{r^{2}}{8}}, \\
\label{332b}&1=-\frac{1}{2}(3\kappa_3-\nu)(1-\frac{r^{2}}{4})e^{-\frac{r^{2}}{8}}.
\end{align}
\end{subequations}
To calculate $r$, \eqref{332a} and \eqref{332b} are divided to each other to have
\begin{equation}\label{333d}
r^{3}-\Delta\mu r^{2}+2Cx_0r^{2}-\frac{4}{3}r-8Cx_0+4\Delta\mu=0.
\end{equation}
$r\approx2+\frac{4}{3\Delta\mu}$ is the solution of \eqref{333d}, of order $O(\frac{1}{\Delta\mu})$. We solve \eqref{332a} for $\kappa_3$ as
\begin{equation}
\kappa_3+\nu=-\frac{r+2Cx_0-\Delta\mu}{r}e^{\frac{r^{2}}{8}}.
\end{equation}
For the value of $r$, this equation is rewritten to find
\begin{equation}
\kappa_3=\frac{1}{2(2+3\Delta\mu)}
         \Big\{-4\nu-6\nu\Delta\mu+\big[3\Delta\mu^{2}-6(1+Cx_0)\Delta\mu-4\big]e^{\frac{1}{8}(2+\frac{4}{3\Delta\mu})^{2}}\Big\}
\end{equation}
\begin{figure}[!ht] \centering
{\label{sek1a} %% label for first subfigure
\includegraphics[scale=1]{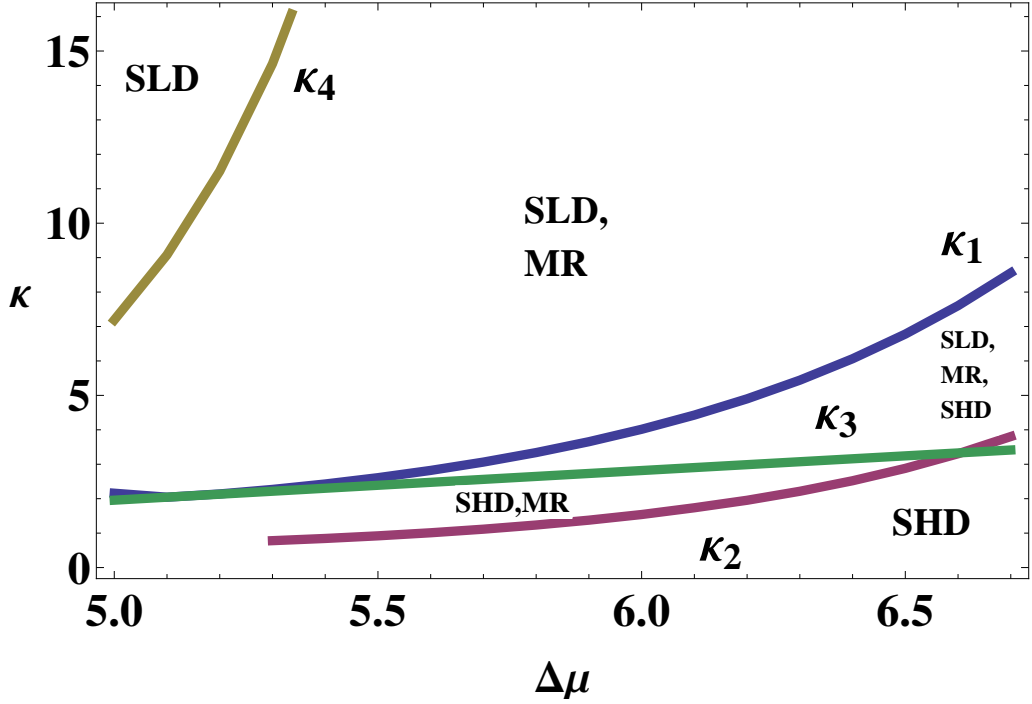}} \hspace{1.5cm}
\caption{{Stability diagram of the triad \eqref{main} with symmetric coupling ($\nu=0$) for chain network}.}
\label{sek1} %% label for entire figure
\end{figure}
In Figure 4, the boundary $\kappa_4$ is the threshold beyond which Majority Rule does not occur. $\kappa_4$ was determined by direct simulation of the dynamical system and detecting the critical values of $\kappa$ versus $\Delta\mu$ that we note in the following table.

\begin{table}[h]
\centering
\caption{Boundary $\kappa_4$}
\label{my-label}
\begin{tabular}{|c|c|c|c|c|c|c|c|c|c|c|}
\hline
$\Delta\mu$&$5$&$5.1$&$5.2$&$5.3$&$5.4$&$5.5$&$5.6$&$5.7$&$5.8$&$5.9$ \\ \hline
$\kappa_4$&$7.19$&$9.08$&$11.51$&$14.63$&$18.65$&$23.84$&$30.58$&$39.33$&$50.76$&$65.72$ \\ \hline
\end{tabular}
\end{table}

\section{Discussions and Future Work}

We want to conclude by pointing out new directions that the analysis in this work can be further taken to.

The system \eqref{sys} is indeed a modification to \eqref{gab}, including the leadership effect. See that \eqref{sys} can be written in the form
\begin{equation}\label{syss}
\frac{dx_i}{dt}=-(\gamma_i+C_i)\Big(x_i-\frac{\gamma_i\mu_i+x_0C_i}{\gamma_i+C_i}\Big)+\sum_{j=1}^N\kappa_{ij}h(x_j-x_i).
\end{equation}
Expressed this way, the leadership effect modifies \eqref{gab} in a way that each agent is under the stress of getting away from the modified bias $\displaystyle \frac{\gamma_i\mu_i+x_0C_i}{\gamma_i+C_i}$, including the opinion of himself and that of the leader. In this picture, we considered small values of $C_i$s so that the profiles SHD, MR and SLD of \cite{Gabbay2014} arise as behavioral patterns of the system and determined how the main results of \cite{Gabbay2014} change in the new case. The bifurcation analysis is done according to the parameters $\Delta \mu$ and $\kappa$. A natural problem would be a bifurcation analysis including the parameter $C_i$. This is a more fundamental question to answer: How the power  $C_i$ of the leader on the agents affect the agreement level of the system? Thus, bifurcation analysis in the parameter spaces $(\Delta \mu,C_i)$, $(\kappa,C_i)$, $(\Delta \mu, \kappa,C_i)$ are possible open problems in this direction.

What are possible behaviours of the system in that case? When the values of $C_i$ and the coupling constants are taken as different then we did here, it is possible to obtain the profiles given in Figures 4 and 5. In Figure 4, all the agents change their opinions a considerable amount compared to their initial opinions; with the final profile SHD in (a), MR in (b) and SLD in (c). The effect of the leader on each agent in the upward direction is quite apparent. This complete control case of \eqref{sys1} is going to be the subject of another paper.
\begin{figure}[!ht] \centering
\subfigure[$\kappa=1$]{\label{sek1a} %% label for first subfigure
\includegraphics[scale=0.50]{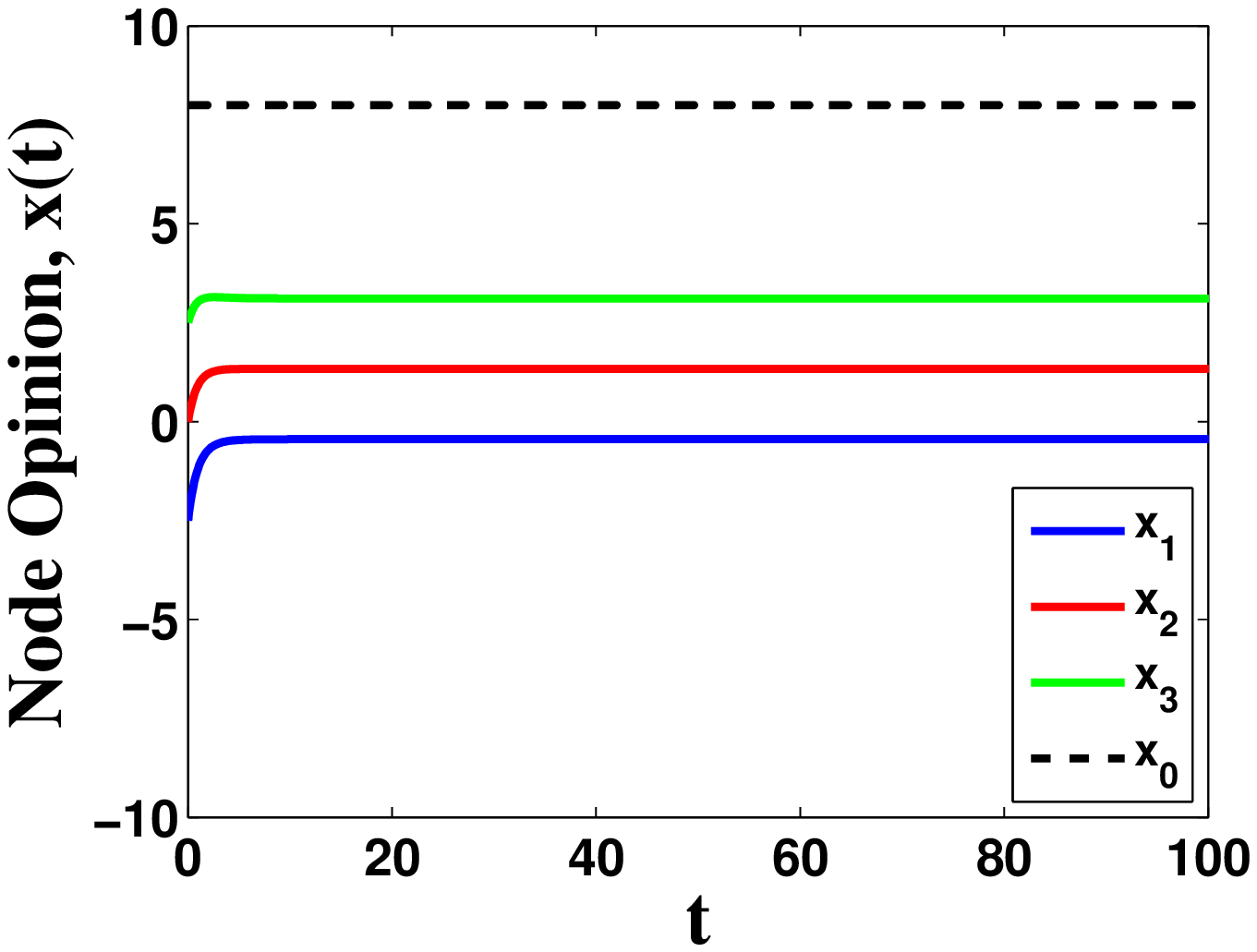}} \hspace{2cm}
\subfigure[$\kappa=1.5$]{\label{sek1b} %% label for second subfigure
\includegraphics[scale=0.50]{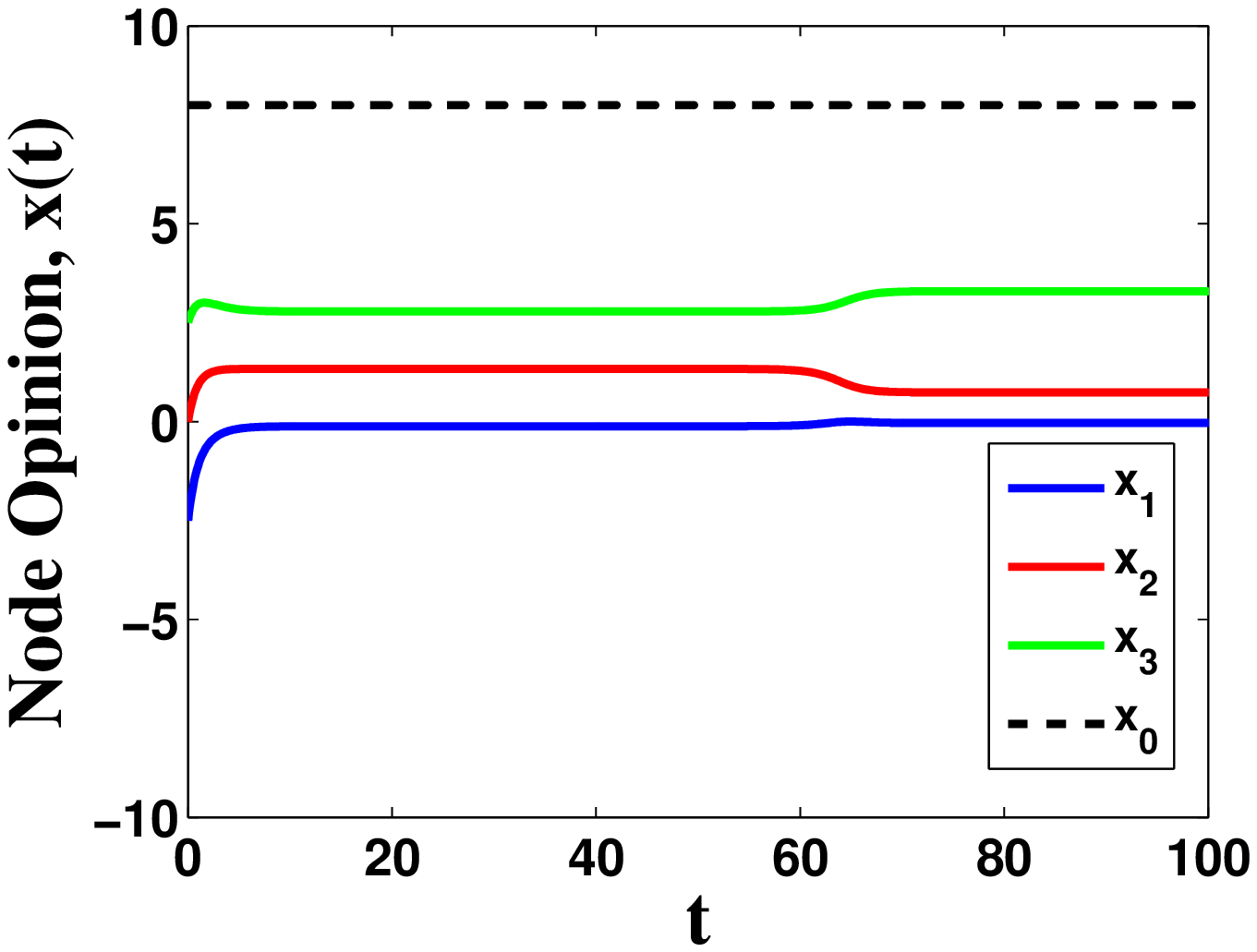}}
\subfigure[$\kappa=4$]{\label{sek1c} %% label for first subfigure
\includegraphics[scale=0.50]{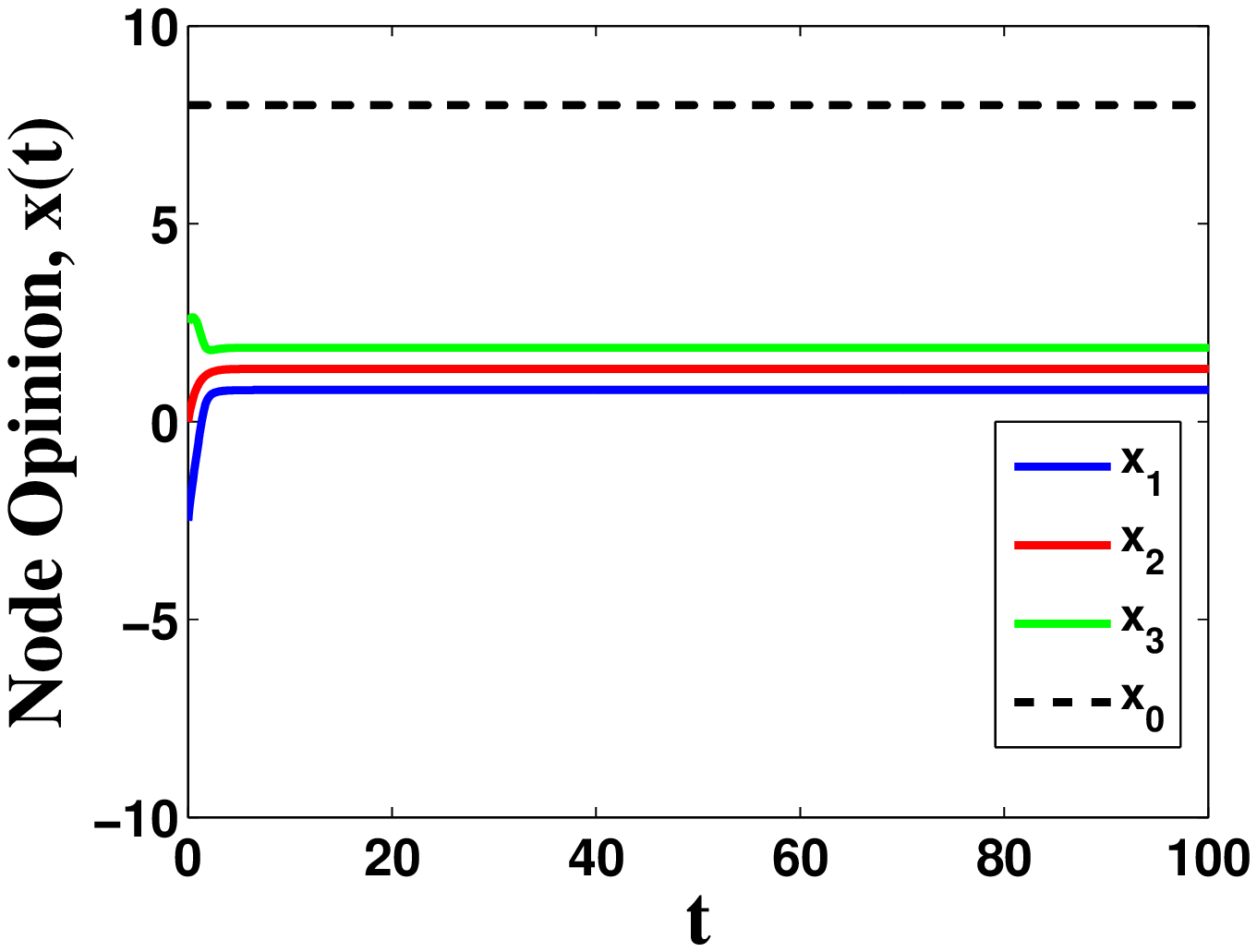}}\hspace{2cm}
 \caption{{Equilibrium outcomes in the dynamical system \eqref{sys1} with complete control case: $C_1=C_2=C_3=C=0.2$, $x_0=8$. Initial disagreement: $\Delta\mu=5$ and initial conditions: $x_1(0)=-2.5$, $x_2(0)=0$, $x_3(0)=2.5$.}}
\label{sek1} %% label for entire figure
\end{figure}

In Figure 5(a), the strength of the leader on the individuals is quick and more dominant, making them reach a common decision very close to the leader's.  Figure 5(b) has been our favorite, as it exhibits the very interesting feature MR first, and an SLD afterwards.
\begin{figure}[!ht] \centering
\subfigure[$\kappa=14$]{\label{sek1a} %% label for first subfigure
\includegraphics[scale=0.45]{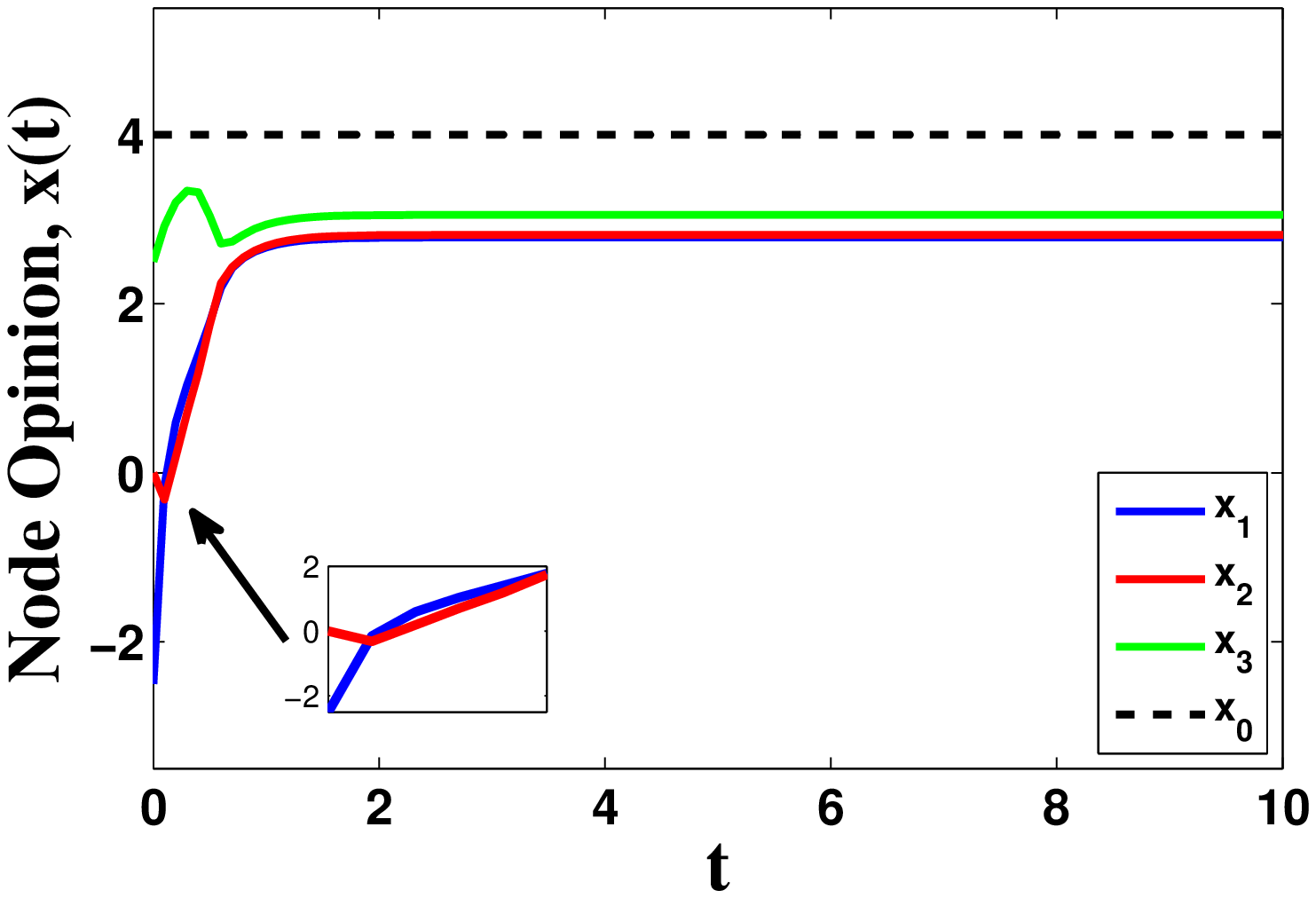}} \hspace{1cm}
\subfigure[$\kappa=1.5$]{\label{sek1b} %% label for second subfigure
\includegraphics[scale=0.45]{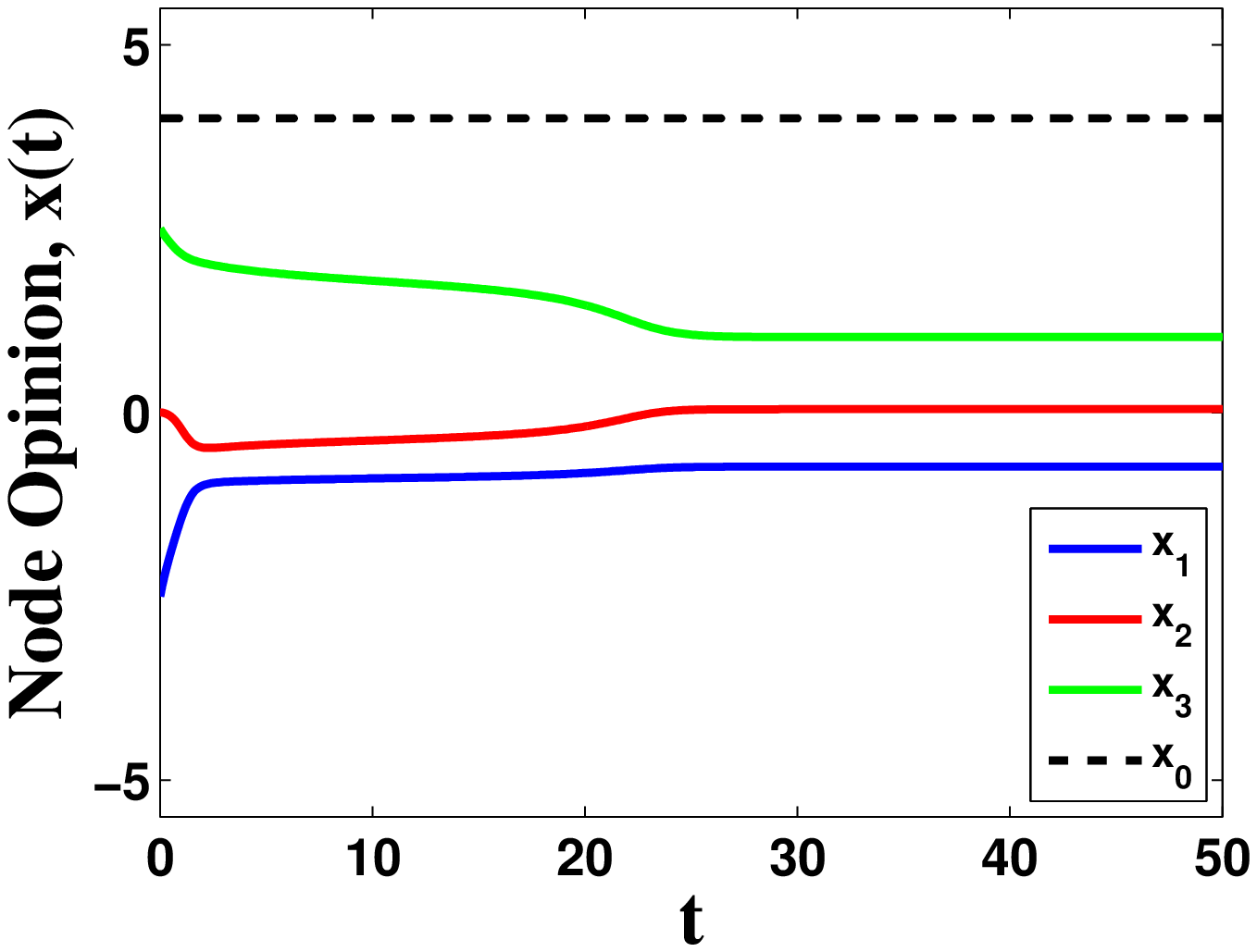}}
 \caption{{Numerical simulation of \eqref{sys1} with end nodes controlled: (a) $C_1=C_3=4$, $C_2=0$ (b) $C_1=-C_3=0.19$, $C_2=0$. Initial disagreement: $\Delta\mu=5$ and initial conditions: $x_1(0)=-2.5$, $x_2(0)=0$, $x_3(0)=2.5$ with $x_0=4$. }}
\label{sek1} %% label for entire figure
\end{figure}

The leadership effect we considered is in the form of a simple harmonic force $C_i(x_0-x_i)$. In literature there are forces of general power types; i.e., we could also consider forces of the type  $C_i|x-x_0|^p$, where the power $p$ can also assume negative values. If $p>0$, the effect of the leader on the agent becomes smaller when the agent's opinion is around the leader's opinion. If $p<0$, the effect of the leader on agents far away are not much; however, if an agent comes around the leader's opinion, the effect is rapidly increasing, due to the reciprocal dependence on the distance. As an example to the literature with this type of effect, we can mention \cite{Kolokolnikov2014}.

\end{document}